\newif\ifHAL
\HALtrue

\documentclass[a4paper,10pt]{article}

\usepackage{a4wide,geometry}
\usepackage{amsmath}
\usepackage{amsmath,amssymb,enumitem,multicol}
\usepackage{mathabx}
\usepackage{amsthm}
\usepackage{upgreek}
\usepackage{bigints}
\usepackage{parskip}
\usepackage{accents}
\usepackage{mathrsfs}
\usepackage{mathtools}
\usepackage{xcolor}
\usepackage{float}
\usepackage[utf8]{inputenc}
\DeclareUnicodeCharacter{2212}{-}
\usepackage[english]{babel}
\usepackage{graphicx}
\usepackage{subcaption}
\usepackage{url}
\usepackage{tikz}
\usepackage{bbm}

\usepackage{amsmath,amsthm,amssymb}

\newtheorem{theorem}{Theorem}[section]

\newtheorem{remark}[theorem]{Remark}

\numberwithin{equation}{section}

\usepackage{accents}

\usepackage{xcolor}
\usepackage{scalerel}

\newcommand{\x}{\hat{x}}

\def\b#1{\boldsymbol{#1}}
\newcommand{\wh}{Warburton--Hesthaven}

\begin{document}

\title{A note on the constants in inverse trace inequalities for polynomials orthogonal to lower-order subspaces}

\author{Zhaonan Dong\footnotemark[1] \footnotemark[2] \and Tanvi Wadhawan\footnotemark[1]  \footnotemark[2]}

\footnotetext[1]{Inria, 48 rue Barrault, 75647 Paris, France}
\footnotetext[2]{CERMICS, CNRS, ENPC, Institut Polytechnique de Paris, 6 \& 8 avenue B.~Pascal, 77455 Marne-la-Vall\'{e}e, France}

\date{}

\maketitle

\begin{abstract}
We derive sharp, explicit constants in inverse trace inequalities for polynomial functions belonging to $\mathbb{P}_p(T)$ (polynomial space with total degree $p$) that are orthogonal to the lower-order subspace $\mathbb{P}_n(T)$, $n\leq p$, where $T$ denotes a $d$-dimensional simplex. The proofs rely on orthogonal polynomial expansions on reference simplices and on a careful analysis of the eigenvalues of the relevant blocks of the face mass matrices, following the arguments developed in~\cite{warburton2003constants}. The novelty is that the extremal face-mass eigenvalue is computed after removing the polynomial modes of degree at most $n$. This yields inverse trace inequality constants involving the factor $(p-n)(p+n+d+1)$ instead of the classical factor $(p+1)(p+d)$, and therefore quantifies the gain in $p$ available in projection-error estimates. These results are very useful in the $hp$-analysis of the hybrid Galerkin methods, e.g. hybridizable  discontinuous Galerkin methods, hybrid high-order methods, etc.

\end{abstract}

\noindent
\textbf{Keywords.} $hp$-analysis, inverse estimate, discrete trace inequality
\medskip\noindent

\textbf{MSC.} 65N30

\section{Introduction}

Inverse inequalities are a cornerstone of finite element theory and play an essential role in the analysis of modern numerical schemes. Classical norm-equivalence arguments, dating back to the early work of Ciarlet~\cite{Ciarlet}, provide general inverse estimates for simplicial finite element spaces. Over the past several decades, there has been substantial progress in deriving sharp and explicit constants in inverse inequalities. We refer, in particular, to the monograph by Schwab~\cite[Section 4.6]{schwab}, where several inverse inequalities are established with $hp$-explicit bounds up to a generic constant that is independent of $h$ and $p$.

For inverse trace inequalities, a significant step forward was taken by Warburton and Hesthaven~\cite{warburton2003constants}, who established constant-free, $h$- and $p$-explicit bounds for inverse inequalities for polynomials in the space $\mathbb{P}_p(T)$, consisting of polynomials of total degree at most $p$ on simplices $T$. In particular, the authors also identified extremal polynomials for which the inverse inequality becomes an identity. These inequalities play a central role in the stability analysis of discontinuous Galerkin methods and other nonconforming methods; see, e.g.,~\cite{DiPietroErnbook,Polydgbook}.

On the other hand,  the results of Warburton and Hesthaven are sharp with respect to the mesh size $h$ and the polynomial degree $p$. However, their results concern the trace of a general polynomial in the space $\mathbb{P}_p(T)$. If one applies their results to a polynomial in $\mathbb{P}_p(T)$ that is orthogonal to lower-order subspaces $\mathbb{P}_n(T)$, the bound can be pessimistic in terms of $p$. This is because the orthogonality with respect to the lower-order subspaces removes the contribution of the lower polynomial modes, which results in a sharper bound in $p$.  Such estimates typically arise in hybrid finite element formulations such as HDG and HHO methods, where traces of local projection errors enter stability estimates and a posteriori error indicators. More precisely, the quantity to be controlled is often of the form $v_T-\Pi^n (v_T)$, or a local lifting generated by face unknowns after subtracting its cell projection. The present estimate yields a constant that reflects the first retained polynomial degree $n+1$, rather than on all polynomial modes ranging from degree zero to degree $p$. As an application, we refer to~\cite{dong:hal-05498158}, where the present result is employed in the derivation of $hp$-\emph{a posteriori} error estimates for hybrid high-order methods.

\paragraph{Contribution of this work.}
In this work, we extend the results from \cite{warburton2003constants} to the setting of inverse trace inequalities for polynomials orthogonal to a lower-order subspace. The main result is the following theorem.
\begin{theorem}[Inverse trace inequality for polynomials orthogonal to lower-order subspaces]\label{Theorem:inverse trace inequality}
Let $T$ be a $d$-dimensional simplex and let $F$ be a $(d-1)$-dimensional simplicial face such that $F \subset \partial T$. The following estimate holds for all $\xi \in \mathbb{P}_p({T})$ with $0\leq n \leq p$:
\begin{align}\label{eq:sharp_estimate}
\|\xi - \Pi^{n}(\xi)\|^2_{L^2({F})}
\leq \frac{(p-n)(p + n+1 + d)}{d} {\frac{|F|}{|T|}}
\|\xi - \Pi^{n}(\xi)\|^2_{L^2(T)},
\end{align}
where ${F}$ denotes a $(d-1)$-dimensional face of $T$ and $\Pi^{n}$ denotes the $L^2$-orthogonal projection onto the $\mathbb{P}_n({T})$.  For $n=-1$, we set $\Pi^{-1}_T (\xi) := 0$.
\end{theorem}

On the other hand, if we apply the classical inverse trace inequality from \cite{warburton2003constants} directly to $\theta := \xi - \Pi^{n}(\xi)$, with $\theta \in \mathbb{P}_p(T)$, we infer
\begin{align*}
\|\theta\|_{L^2(F)}^2 \le \frac{(p+1)(p+d)}{d}\frac{|F|}{|T|}\|\theta\|_{L^2(T)}^2.
\end{align*}
However, $\theta$ is orthogonal to $\mathbb{P}_{n}(T)$, which removes the contribution of the low-order polynomial modes. Consequently, the analysis provides the sharper bound
\begin{align*}
\|\theta\|_{L^2(F)}^2 \le \frac{(p-n)(p + n+1 + d)}{d}\frac{|F|}{|T|}\|\theta\|_{L^2(T)}^2.
\end{align*}
It is easy to see that the new inverse trace inequality bound is always strictly sharper than the bound using \cite{warburton2003constants} for $0\leq n \leq p$. For the case $n=-1$, the new bound is equal to the classical bound.
\par
\noindent
Compared with~\cite{warburton2003constants}, the analysis below preserves the original orthogonal basis and instead exploits the orthogonality condition $\theta \perp \mathbb{P}_n(T)$ through a spectral characterization of the associated truncated coefficient space. The main observation is that the truncation shifts the relevant extremal eigenvalue from a contribution involving sum over all polynomial degrees $0,\ldots,p$ to a sum over the retained degrees $n+1,\ldots,p$. For clarity, we first present the argument in  two dimension before establishing the general result on a $d$-simplex in Theorem~\ref{Theorem:inverse trace inequality}.

The remainder of the paper is organized as follows. Section~\ref{sec:notation} introduces notation and preliminary concepts. Section~\ref{sec:2D} present the analysis on the  triangles. In Section \ref{sec: d-dimensional simplex}, we generalize the result to the $d$-dimensional simplices. Section~\ref{sec:numerics} gives illustrative numerical comparisons of the constants.

\section{Preliminaries and notation}\label{sec:notation}
Let us start by introducing  the following notations used throughout the paper:
\begin{itemize}
\item Let  $T$ denote a $d$-dimensional simplex with boundary $\partial T$, and  let $F$ be a $(d-1)$-dimensional simplicial face such that $F \subset \partial T$. Equivalently, $T$ is the affine image of the reference simplex
\[
\hat{T}_d:=\left\{(r_1,r_2,\dots,r_d)\;:\; |r_i|\leq 1, \sum_{i=1}^d r_i \leq 2-d\right\}.
\]
\item  For $p \ge 0$, $\mathbb{P}_p(T)$ is the space of polynomials of total degree at most $p$ on $T$. 
\item  Any $u \in \mathbb{P}_p(T)$ can be expanded in an $L^2(T)$-orthonormal basis $\{\psi_{\b{i}}\}_{|\b{i}|=0}^{p}$:
\begin{align*}
u(\b{x}) = \sum_{|\b{i}|=0}^{p}  u_{\b{i}} \, \psi_{\b{i}}(\b{x}),
\qquad  u_{\b{i}} = \int_T u(\b{x}) \psi_{\b{i}}(\b{x}) \, d\b{x},
\end{align*}
where $\b{x}:=(x_1,x_2,\dots, x_d)^{\bot}$.
\item  The $L^2$ norms on elements and faces are denoted by $\|\cdot\|_{L^2(T)}$ and $\|\cdot\|_{L^2(F)}$, respectively.
\item For an integer $0\leq n\leq p$, let $\Pi^{n}$ denote the $L^2$-projection onto $\mathbb{P}_{n}(T)$. For $n=-1$,  we set $\Pi^{-1} (\xi) := 0$.
\end{itemize}

\section{Inverse inequality on a 2D triangle}\label{sec:2D}
In this section,  we will derive the explicit form of the constant in Theorem \ref{Theorem:inverse trace inequality} on a 2-dimensional simplex.

Consider the reference triangle
$\hat{T}_2 := \{(r_1, r_2) \mid -1 \leq r_1,r_2 \leq 1,\ r_1+r_2 \leq 0\}$.
To parametrize $\hat{T}_2$, we introduce the Duffy mapping with coordinates $a,b \in [-1,1]$ via
\begin{align*}
r_1 &= \frac{(1+a)}{2}(1-b) - 1, \quad
r_2 = b, \quad -1 \leq a,b \leq 1,
\end{align*}
which maps the unit square $[-1,1]^2$ onto $\hat{T}_2$.
We employ the orthonormal polynomial basis on $\hat{T}_2$ introduced by Proriol~\cite{proriol57} and subsequently refined by Koornwinder~\cite{koornwinder75} and Dubiner~\cite{dubiner91},
indexed by integer pairs $(i,j)$ with $i+j \leq p$:
\begin{align*}
\psi_{ij}(r_1, r_2) =
\bigg(\frac{P_i^{(0,0)}(a)}{\sqrt{\frac{2}{2i+1}}}\bigg)
\bigg(\frac{\left(\frac{1-b}{2}\right)^i P_j^{(2i+1,0)}(b)}{\sqrt{\frac{1}{i+j+1}}}\bigg),
\end{align*}
where $P_n^{(\alpha,\beta)}(x)$ denotes the $n$-th order Jacobi polynomial on $[-1,1]$.
Any polynomial $\xi \in \mathbb{P}_p(\hat{T}_2)$ admits the expansion
\begin{align*}
\xi(r_1, r_2) = \sum_{i+j \leq p} \xi_{(i,j)} \psi_{ij}(r_1, r_2).
\end{align*}
Since $\theta = \xi - \Pi^{n} (\xi)$ is orthogonal to $\mathbb{P}_n(\hat{T}_2)$, it admits the expansion
\begin{align*}
\theta(r_1, r_2) = \sum_{i+j \geq n+1}^{p} \xi_{(i,j)} \psi_{ij}(r_1, r_2).
\end{align*}

Focusing on the face $\hat{F}$ where $r_2=-1$ (equivalently $b=-1$), we compute
\begin{align}
\int_{\hat{F}} \theta^2(r_1,-1)\, dr
&= \int_{-1}^{1} \theta^2(a,-1)\, da
= \Theta^{\top} L \Theta, \label{eq:2D-edge-integral}
\end{align}
where $\Theta$ is the vector of coefficients $\xi_{(i,j)}$ with $i+j \in [n+1,p]$,
and $L$ is the associated face mass matrix with entries
\begin{align*}
L_{(ij)(kl)}
&= \int_{-1}^{1} \psi_{ij}(r_1,-1) \psi_{kl}(r_1,-1)\, dr \\
&= \delta_{ik} \, (-1)^{j+l} \,
\sqrt{i+j+1} \,
\sqrt{k+l+1},
\end{align*}
where $\delta_{ik}$ follows from the $L^2$-orthogonality of the Legendre polynomials $P_i^{(0,0)}(\x)$.

The matrix $L$ is block diagonal, with blocks indexed by $i=0,\ldots,p$. Hence, its spectral radius can be determined by analyzing each block separately. For the $i$-th block, $j$ ranges from $\max\{0,n+1-i\}$ to $p-i$, each block is a rank-one matrix of the form:
\begin{align*}
Z^{(i)} = \mathbf{v}^{(i)} (\mathbf{v}^{(i)})^\top, \quad
\mathbf{v}^{(i)}_j = (-1)^j \sqrt{i+j+1}.
\end{align*}
The spectral radius $\rho(L)$ of $L$ equals the maximum of the eigenvalues across all blocks. For $0\le i\le n+1$, the eigenvalue of the $i$-th block $Z^{(i)}$ is
$$
\lambda^{(i)} = \sum_{j=n+1-i}^{p-i} (i+j+1) =  \frac{(p-n)(p+n+3)}{2},
$$
which is independent of $i$. For $i>n+1$, the lower index is $j=0$, and
\[
\lambda^{(i)}
=\sum_{j=0}^{p-i}(i+j+1)
\le \frac{(p-n)(p+n+3)}{2}.
\]
Indeed, writing $i=n+1+q$ with $q\ge0$, one verifies that the difference between the right-hand side and $\lambda^{(i)}$ is
$\frac12 q(2n+3+q)\ge0$. Consequently, the spectral radius $\rho(L)$ is
\begin{align*}
\rho(L) =   \frac{(p-n)(p+n+3)}{2}.
\end{align*}
Using \eqref{eq:2D-edge-integral}, we obtain
\begin{align*}
\|\theta\|_{L^2(\hat{F})}^2 \leq \rho(L) \|\theta\|_{L^2(\hat{T}_2)}^2
=  \frac{(p-n)(p+n+3)}{2} \|\theta\|_{L^2(\hat{T}_2)}^2.
\end{align*}
Finally, a standard scaling argument yields \eqref{eq:sharp_estimate} for $d=2$.

\section{Inverse inequality on the $d$-dimensional simplex}\label{sec: d-dimensional simplex}
In this section, we establish the inverse trace inequality on a general $d$-dimensional simplex, thereby proving the main result of the article, namely Theorem~\ref{Theorem:inverse trace inequality}.

Consider the reference $d$-simplex $\hat{T}_d:=\left\{(r_1,r_2,\dots,r_d)\;:\; |r_i|\leq 1, \sum_{i=1}^d r_i \leq 2-d\right\}$. Next, to parametrize $\hat{T}_d$, we introduce the Duffy mapping with coordinates $a_i \in [-1,1]$, for $i=1,\dots,d$,  such that
\begin{equation*}\begin{split}
r_1 &= \frac{(1+a_1)}{2}  \frac{(1-a_2)}{2} \frac{(1-a_3)}{2} \dots (1-a_d)-1,  \\
r_2 &= \frac{(1+a_2)}{2}  \frac{(1-a_3)}{2} \frac{(1-a_4)}{2} \dots (1-a_d)-1,  \\
\vdots & ,\\
r_d &= a_d,
\end{split}\end{equation*}
which maps the the bi-unit cube $[-1,1]^d$ onto the reference simplex $\hat{T}_d$.
We employ the orthonormal polynomial basis for $\hat{T}_d$ indexed by multi-indices $\b{i}=(i_1,i_2,\dots, i_d)$ as introduced in \cite{warburton2003constants} as
\begin{equation*}\begin{split}
\psi_{(\b{i})}(r_1,r_2,\dots,r_d) &= \bigg(\frac{P_{i_1}^{(0,0)}(a_1)}{\sqrt{\frac{2}{2i_1+1}}}\bigg) \prod_{l=2}^{d}\bigg(\frac{\left(\frac{1-a_l}{2}\right)^{N_l(\b{i})}
P_{i_l}^{(2N_l(\b{i})+l,0)}(a_l)}{{\sqrt{\frac{2}{2N_l(\b{i})+d}}}}\bigg) ,
\end{split}\end{equation*}
where $N_l(\b{i}) = \sum_{j=1}^{l}i_j$ and $P_n^{(\alpha,\beta)}(x)$ denotes the \(n\)-th order Jacobi polynomial defined on \([-1,1]\). Any polynomial $\xi \in \mathbb{P}_p(\hat{T}_d)$ admits the expansion
$\xi(r_1,r_2,\dots,r_d) = \sum\limits_{|\b{i}|=0}^p {\xi}_{(\b{i})} \, \psi_{\b{i}}(r_1,r_2,\dots,r_d)$. Since $\theta=\xi-\Pi^n (\xi)$ is orthogonal to $\mathbb P_n(\hat{T}_d)$, it admits the expansion
\begin{equation*}\begin{split}
\theta(r_1,r_2,\dots,r_d) = \sum_{|\b{i}|= n+1}^{p} {\xi}_{( \b{i})} \psi_{\b{i}}(r_1,r_2,\dots,r_d).
\end{split}\end{equation*}
Next, we focus on the face $\hat{F}$ where $a_d=-1,$ and evaluate the corresponding face integral:
\begin{align}\label{eqn2:Lemma2.3}
\int\limits_{\hat{F}} \theta^2(r_1,r_2,\dots,-1)~d{r_1}d{r_{2}}~d{r_{d-1}} = \int\limits_{a_d=-1} \theta^2(a_1,a_2,\dots,-1)~d{a_1}d{a_{2}}~d{a_{d-1}} \;=\; \Theta^{\top} L \Theta,
\end{align}
where $\Theta$ is the vector of coefficients ${\xi}_{(\b{i})}$ corresponding to $|\b{i}| \in [n+1,p]$, and $L$ denotes the associated face mass matrix. The entries of the face mass matrix $L$ are given by
\begin{equation*}\begin{split}
L_{(\b{i})(\b{j})}
=\prod_{l=1}^{d-1} \delta_{i_l j_l}
(-1)^{(i_d + j_d)}\,
\sqrt{\frac{2N_d(\b{i}) +d}{2}}
\;\sqrt{\frac{2N_d(\b{j}) +d}{2}} \,.
\end{split}\end{equation*}
Consequently, $L$ is block diagonal with respect to the first $d-1$ indices. For a fixed multi-index $(i_1,\ldots,i_{d-1})$, set $\alpha=\sum_{l=1}^{d-1}i_l.$
Then, within the corresponding block, $k=i_d$ ranges from $\max\{0,n+1-\alpha\}
\quad\text{to}~
p-\alpha.$
Each block is rank one and can be written in the form
\begin{equation*}
Z_{(\alpha)}= {\b{\nu}}^{(\alpha)}({\b{\nu}}^{(\alpha)})^{\top}, ~\text{with} \quad {\b{\nu}}^{(\alpha)}_k = (-1)^k \sqrt{\frac{2(\alpha+k)+d}{2}}.
\end{equation*}

The maximum eigenvalue of each rank-one block equals the squared norm of \(\b{\nu}^{(\alpha)}\):
\begin{equation*}\begin{split}
\lambda^{(\alpha)} = \sum_{k=n+1-\alpha}^{p-\alpha} \frac{2(\alpha+k)+d}{2}=  \frac{(p-n)(p+n+d+1)}{2},
\end{split}\end{equation*}
for $0\le \alpha\le n+1$. If $\alpha>n+1$, the sum starts from $k=0$, and the resulting value is no larger than the one displayed above. Therefore, the spectral radius of $L$ is attained by the blocks corresponding to $0\le \alpha\le n+1$, and hence
$$
\rho(L) =   \frac{(p-n)(p+n+d+1)}{2}.
$$
If $n=-1$, no modes are removed and the estimate reduces to the classical
Warburton--Hesthaven constant. Combining \eqref{eqn2:Lemma2.3} with the above spectral estimate yields
\begin{equation*}\begin{split}
\|\theta\|_{L^2(\hat{F})}^2 \leq \rho(L) \|\theta\|_{L^2(\hat{T}_d)}^2 = \frac{(p-n)(p+n+d+1)}{2}\|\theta\|_{L^2(\hat{T}_d)}^2.
\end{split}\end{equation*}
Finally, a standard scaling argument yields the desired estimate \eqref{eq:sharp_estimate}.
\begin{remark}
For $n=-1$, no polynomial modes are removed, and the resulting constant coincides with the classical Warburton--Hesthaven constant on the $d$-dimensional simplices.
\end{remark}

\section{Numerical examples}\label{sec:numerics}
In this section,
we provide a simple numerical illustration of the improvement in the constants on the reference triangle $\hat{T}_2$. Let
\[
C_{\rm new}(p,n,2):=\frac{(p-n)(p+n+3)}{2}, \qquad
C_{\rm WH}(p,2):=\frac{(p+1)(p+2)}{2} = C_{\rm new}(p,-1,2).
\]
The quantity $C_{\rm WH}(p,2)$ corresponds to the constant obtained by applying the \wh{} inverse trace inequality directly to functions in $\mathbb{P}_p(\hat{T}_2)$, whereas $C_{\rm new}(p,n,2)$ exploits the additional orthogonality condition $\theta\perp\mathbb{P}_n(\hat{T}_2)$.

We compute the largest eigenvalue $\rho$ of the face mass matrix using the Matlab function \texttt{eigs} for the cases $n=0$ and $n=p-1$, with $p\in\{0,\ldots,30\}$. The results are reported in Table~\ref{tab:constants}. We observe that, for $n=0$, the computed eigenvalues coincide with

\[
C_{\rm new}(p,-1,2)= C_{\rm WH}(p,2)=\frac{(p+1)(p+2)}{2}.
\]
Moreover, for $n=p-1$, the computed eigenvalues coincide with
\[
C_{\rm new}(p,p-1,2)=p+1,
\]
thereby confirming the theoretical result.
\begin{table}[htbp]
\centering
\caption{Values of $C_{\rm WH}(p,2)$ and $C_{\rm new}(p,p-1,2)$.}
\label{tab:constants}
\begin{tabular}{|c|c|c|c|c|c|}
\hline
$p$ & $C_{\rm WH}(p,2)$ & $C_{\rm new}(p,p-1,2)$ &
$p$ & $C_{\rm WH}(p,2)$ & $C_{\rm new}(p,p-1,2)$ \\
\hline
1  & 3   & 2  & 16 & 153 & 17 \\
2  & 6   & 3  & 17 & 171 & 18 \\
3  & 10  & 4  & 18 & 190 & 19 \\
4  & 15  & 5  & 19 & 210 & 20 \\
5  & 21  & 6  & 20 & 231 & 21 \\
6  & 28  & 7  & 21 & 253 & 22 \\
7  & 36  & 8  & 22 & 276 & 23 \\
8  & 45  & 9  & 23 & 300 & 24 \\
9  & 55  & 10 & 24 & 325 & 25 \\
10 & 66  & 11 & 25 & 351 & 26 \\
11 & 78  & 12 & 26 & 378 & 27 \\
12 & 91  & 13 & 27 & 406 & 28 \\
13 & 105 & 14 & 28 & 435 & 29 \\
14 & 120 & 15 & 29 & 465 & 30 \\
15 & 136 & 16 & 30 & 496 & 31 \\
\hline
\end{tabular}
\end{table}

\bibliographystyle{siam}
\bibliography{biblio}

\end{document}